\documentclass{amsart}
\usepackage[T1]{fontenc}
\usepackage{amssymb}
\usepackage[hyperref]{degt}

\def\2{\color{red}}
\def\3{\color{green}}
\def\4{\color{blue}}

\def\Dg:{{\bf Dg:\enspace}\ignorespaces}

\def\KK-{\pdfstr{K3}{$K3$}-\penalty0\ignorespaces}

\def\0{}

\def\dual{^\vee}

\let\graph\Gamma

\let\Weyl\Delta
\def\fp{\Weyl^{\sharp}}

\def\Fano{\Cal{F}}
\def\Fano{\frak{F}}

\def\Fn{\operatorname{Fn}}

\def\root{\operatorname{root}}

\def\FfNS{\Fano_{\ttD}^d(\graph)}

\def\Qano{\frak{G}}
\def\barFfNS{\Qano_{\ttD}^d(\graph)   }

\def\Gapp{\omega}

\def\cl[#1]{[\![#1]\!]}

\def\NS{\operatorname{NS}}

\def\qgen{spanned }

\def\resp.{resp\PERIOD}
\let\geq\ge
\let\leq\le

\let\M=M
\let\o=o

\let\ex=e
\let\iso=p

\def\basis#1{\{#1\}}

\def\CK{\Cal K}

\def\bL{\bold{L}}

\def\CC{\Cal C}

\def\bL{\bold L}
\def\CU{\Cal U}
\def\CV{\Cal V}
\def\CK{\Cal K}
\def\CN{\Cal N}

\def\Br{\operatorname{Br}}
\def\dual{^\vee}
\let\kk\varkappa

\def\ba{\bar{a}}
\def\sa{a^*}
\def\units{^{\times\!}}

\newcommand{\rr}{\mathfrak r}
\newcommand{\ttH}{h}  
\newcommand{\ttD}{2D} 
\def\ttD{{2\ooD}}
\def\ooD{n}

\def\CD{\Cal{D}}
\def\CM{\Cal{M}}
\def\CK{\Cal{K}}

\def\ttA#1{\pdfstr{AA}{\tA}\sb{#1}}

\def\p#1{#1^+}
\def\ppart#1{\operatorname{Ln}_{#1}}
\def\ppart#1{\nu_{#1}}

\newcounter{row}
\def\therow{\arabic{row}}

\def\vtab{\vbox\bgroup\def\(##1){[##1]}%
 \ialign\bgroup\hss\refstepcounter{row}\ \therow##\hss\quad&
  \strut\hss$##$\hss\quad&$\singset{##}$\ \hss\cr
  \noalign{\hrule\vspace{2pt}}
  \omit\ \hss$\#$\hss\quad&\text{ord}&\omit\text{pattern}\hss\cr
  \noalign{\vspace{1pt}\hrule\vspace{2pt}}
  }
\def\endvtab{\crcr\noalign{\vspace{1pt}\hrule}\egroup\egroup}

\newcommand{\same}{\simeq}
\newcommand{\sameg}{\simeq}

\title[Periodic sharpness of Miyaoka's bound]
{ Periodic sharpness of Miyaoka's bound for smooth rational curves}

\author{Alex Degtyarev}

\address{%
Department of Mathematics\\
Bilkent University\\
06800 Ankara, TURKEY}

\email{
degt@fen.bilkent.edu.tr}

\thanks{%
A.D. was partially supported by the T\"{U}B\DOTaccent{I}TAK grant 123F111%
}

\author{S\l awomir Rams}
\address{Theoretical Computer Science Department, Faculty of Mathematics and Computer Science, Jagiellonian University,
	ul.\ {\L}ojasiewicza 6,  30-348 Krak\'{o}w, Poland}
\email{slawomir.rams@uj.edu.pl}

\thanks{%
    S.R. was supported by the National Science Centre, Poland, Opus  grant
	no.\ 2024/\allowbreak53/\allowbreak B/\allowbreak ST1/\allowbreak01413 and
partially supported by
{T\"{U}B\DOTaccent{I}TAK B\DOTaccent{I}DEB
2221 Visiting Scientist Fellowship Program}%
}



\begin{document}
\newcommand{\smooth}{}
\begin{abstract}
	We determine the maximal number of smooth
	rational degree~$d$ curves on a \smooth complex $K3$-surface of degree $\ttD$ provided $\ooD$ is sufficiently large as compared to~$d>1$
	(\autoref{thm.main}). We obtain precise characterization of configurations of rational degree~$d$ curves
	 for which Miyaoka's bound is sharp for $nd$  odd and
	 $\ooD \gg d$
	 (\autoref{cor-miyaoka-classification}).
\end{abstract}

\maketitle

\section{Introduction}\label{S.intro}

Configurations  of lines and conics
 on various classes of surfaces
have been studied for over 170 years, with numerous partial results spread over the literature. Within the last decade our understanding
 of line configurations on polarized \KK-surfaces over
various fields became almost complete (\cf.
\cite{degt:lines,degt:sextics,DIS,Degtyarev.Rams.triangular,rams.schuett,rams.schuett:char3}
and further references therein)
and
much has been shown for conics (see
\cite{Degtyarev22conics800,
	Degtyarev23conicsKummer,Degtyarev21conicssextics,
	Degtyarev24conicsBB,degtyarev24conicsquartics,
Miyaoka:lines,	rams.schuett.smoothodd}).

Quite unexpectedly, one can use the Miyaoka--Yau-=Sakai inequality \cite{Miyaoka:orbibundle} to give a very concise description of the behaviour
of the maximal number of degree~$d$ (irreducible, but not necessarily smooth) rational curves on high-degree complex \KK-surfaces.
Indeed, for fixed
integers  $d \geq 1$ and $\ooD \geq 2$ let us define
$$
\rr(\ttD,d) := \max \bigl\{ r_d(X) \bigm|  \mbox{$X$ is a
	degree~$\ttD$ \smooth complex $K3$-surface} \bigr\},
$$
where
$r_d(X)$
stands for the number of degree~$d$ (ireducible, but not necessarily smooth) rational curves on the  surface $X$.
Then, \cite[Proposition A]{Miyaoka:lines})
implies that
\begin{equation} \label{eq-miyaoka}
	\rr(\ttD,d)  \leq 24 \quad\mbox{for}\quad \ooD > 50 d^2 \, .
\end{equation}
If we assume  $\ooD$ sufficiently large  as compared to~$d$,
these
bounds  are sharp for $d \geq 3$
(see \cite[Theorem~1.1]{rams.schuett:24}) and $d=2$ (see \cite[Corollary~1.2]{rams.schuett.smoothodd}), whereas
\cite[Theorem~1.5]{degt:lines} yields
\begin{equation} \label{eq.max.lines}
	\rr(\ttD,1) \in \{21, 22, 24\} \quad \mbox{for } \ooD \gg 0.
\end{equation}
(Precise arithmetic conditions  on the degree $\ttD$ of the polarization
for Miyaoka's bound of 24 lines to be attained can be found in \cite[Theorem~1.5]{degt:lines}.) Consequently,
\eqref{eq-miyaoka} fails to be sharp  for infinitely many $\ooD$ only for lines.

Until recently, far less has been  known
about maximal configurations of \emph{smooth}  rational curves of
small degrees (even for $d=3$, \ie, for twisted cubics).
The
classical construction of \KK-surfaces with many degree~$d$ rational
curves in \cite[Theorem~1.1]{rams.schuett:24} always results in
configurations of {\4 irreducible} rational curves with singularities. In contrast, the recent
construction from \cite{rams.schuett.smoothodd} yields configurations of
smooth curves, but it works only when either $\ooD$ or $d$ is even (and
in the extra
case $2n=2d^2\bmod6$ {\4- see \cite[Proposition~6.1]{rams.schuett.smoothodd}).}

In the present paper
we complete the picture of
the asymptotic behaviour of the maximal
number of \emph{smooth} degree~$d$ rational curves on \smooth
\KK-surfaces of high degree. More precisely, we prove the following analogue of
  \cite[Theorem~1.5]{degt:lines}.
In the statement, we need a very strange arithmetic condition;
for a prime~$q$ and integer~$m$, we denote by
$\ppart{q}(m):=\max\bigl\{r\in\N\bigm|q^r\mathbin|m\bigr\}$
the \emph{$q$-adic valuation} of~$m$.
As usual, this extends to~$\Q$ \via\ $\ppart{q}(m/n):=\ppart{q}(m)-\ppart{q}(n)$.

\theorem[see \autoref{proof.main}] \label{thm.main}
Let $d>1$ be a fixed integer.
If $\ooD\gg d$ {\rm(}\ie~$\ooD\gg0$ as compared to~$d${\rm)}
then the maximal number of smooth
rational degree~$d$ curves on a \smooth complex $K3$-surface of degree $\ttD$ equals
\roster*
\item
$24$, if
either $\ttD d=0\bmod4$ or $\ttD=2d^2\bmod3$\rom;
\item
otherwise, $22$ if $\ttD=4d^2\bmod5$
and, for each prime $q=3\bmod4$, either
$\ppart{q}(\ooD)\ge\ppart{q}(d)$ or $\ppart{q}(\ttD)$
is even\rom;
\item
otherwise, $21$.
\endroster
\endtheorem

\remark\label{rem.period}
The
periodicity of the maximal number of smooth degree~$d$ rational curves
depends on~$d$. For example, of $d$ is square free, the period is
$30=2\cdot3\cdot5$. In general, it is
\[*
30\prod q^{2\left\lfloor\ppart{q}(d)/2\right\rfloor},\quad q=3\bmod4.
\]
If $9\mathrel|d$, the coefficient~$30$ can be reduced down to~$10$.
\endremark

As a
{consequence,}
we obtain
{the}
classification of the maximal
configurations of irreducible rational degree~$d$ curves on \KK-surfaces of
high degree $2\ooD$ with both $\ooD$ and $d$ odd, which in this case
provides a complete answer to the question posed in
\cite[$\S$~3]{Miyaoka:lines} (see  \autoref{rem:even.nd} below for the
discussion of the case of $nd$ \emph{even}).

\corollary[see \autoref{proof.cor-miyaoka-classification}]  \label{cor-miyaoka-classification}
Let $\ooD\gg d$, with  $\ooD d$ odd, and let
 a \smooth complex \KK-surface of degree $\ttD$
 contain
$24$
 irreducible
 rational degree~$d$ curves.
 If   $n \neq  d^2 \bmod 3$, then the curves in question form
 the configuration $24\tA_0^*$.  Otherwise,
  their  configuration is either $8\tA_2$
{\rom(and then the curves are smooth\rom)}
or $24\tA_0^*$.
  \endcorollary

Observe {also}
that a \smooth complex \KK-surface of degree $\ttD$
with  the configurations  given above  exists for any pair of integers $(n,d)$ that satisfy the assumptions of \autoref{cor-miyaoka-classification} (see \cite[$\S$~10]{rams.schuett:24} and \autoref{th.8A2} below).

The arguments
used in \cite{degt:lines}  (and their generalizations in
\cite{rams.schuett:24}) suggest that the maximal number of smooth rational
degree~$d$ curves on a \smooth complex \KK-surface od degree $\ttD$ should be
a periodic function of $\ooD$, with the  period  proportional to the degree
$d$ of the curves in question (see Proposition~\ref{prop:periodic} below), so
its precise description for all values of $d$ may be out of reach. In
contrast to this expectation,
Theorem~\ref{thm.main} shows that the behaviour of the maximal configurations
for most
odd values of $d$
is very much similar to the case $d=1$ of lines.

At early stages of this project our goal was understanding the behaviour
of the
maximal configurations
realizing Miyaoka's bound \cite[Proposition A]{Miyaoka:lines}
for small odd integers~$d$. There are four
major reasons
 why one is able to prove an analogue of  \eqref{eq.max.lines} for \emph{every} odd $d$.
Firstly, by \cite{degt:lines} (see also \cite[Corollary~6.4]{rams.schuett:24}),
for
$\ooD$ large,
we can restrict our attention to
the configurations of rational curves that correspond to
the so-called parabolic graphs. Secondly, for each parabolic graph of interest we
are able to find a basis of the corresponding lattice
(sometimes upon tensoring the latter by~$\Z_p$)
that enables us to control
its discriminant forms
for all
pairs $(\ooD,d)$. Thirdly, the fact that all exceptional
classes in the discriminant group must be orthogonal to the
classes of \emph{both} the
polarization and the fiber (see \autoref{lem.k.perp})
lets us exclude most finite index extensions of the lattices associated
to the graphs
considered
(\cf. \autoref{s.12A1}).  Last but not least, the assumption that $d$ and $\ooD$ are odd
substantially simplifies the structure of
the $2$-torsion of the discriminant.

The paper is organized as follows.
In \autoref{S.prelim}
we recall basic facts
used
and introduce the notation.
We illustrate this by
a very detailed study {\4 of} the lattice associated  to the parabolic graph $\graph := 12\tA_1$ (see \autoref{s.12A1}).
Then,
\autoref{S.A1}
through
\autoref{S.4A4+2A1} are devoted to
finding the degrees in which
various parabolic graphs are geometric.
Finally, in \autoref{S.proof} we give the proof of \autoref{thm.main}.

\remark \label{rem:basefield}
One can ask whether
\autoref{thm.main} holds for \KK-surfaces defined over any algebraically
closed field of characteristic $p \geq 5$ (\cf. \cite[$\S$~11]{rams.schuett:24}). We
presume that some extra cases may appear for small $p>0$ but,
since
our main inspiration was Miyaoka's paper \cite{Miyaoka:lines} that
deals with \emph{complex} surfaces,
we
confine ourselves to the complex case.
\endremark



Although \GAP~\cite{GAP4} was a useful tool helping us guess the statements
and verify our findings, we managed to present the ultimate proofs in a human
readable form, without a reference to software.

\subsection{Acknowledgements}
{This paper
was conceived during Alex Degtyarev's
visit to the Jagiellonian
University (Cracow, Poland). A.D.\ is profoundly grateful to this institution for the
hospitality and excellent working conditions.}
Substantial part of
{the} paper was written during S{\l}awomir Rams' stay at
the Department of Mathematics of Bilkent University (Ankara, Turkey) within
the scope of T\"UB\DOTaccent ITAK ``2221-Fellowships for visiting scientists'' programme.
S.R.\ would like to thank Bilkent University and T\"UB\DOTaccent ITAK for creating perfect
working conditions and numerous inspiring discussions with members of the
Department of Mathematics of Bilkent University.
Both authors are would like to thank Matthias Sch\"utt for a number of
inspiring discussions.


\section{Preliminaries} \label{S.prelim}

In this note we work over the field $\C$ of complex numbers.
Whenever we speak of a \KK-surface of degree $\ttD$,
{or a ($2n$-)polarized $K3$-surface,}
we mean a pair $(X,h)$, where $X$ is a projective
{simply} connected complex
surface with trivial canonical bundle and $h \in \NS(X)$
is \emph{very ample} with $h^2=\ttD$. Sometimes, for the sake of
brevity,
we write ``a \KK-surface'' when we mean ``a \emph{polarized}  \KK-surface''.

We freely extend modular arithmetic to fractions with the denominator prime to
the modulus, so that, \eg, $1/2=2\bmod3$ and $1/2=3\bmod5$, as in \GAP~\cite{GAP4}.


\subsection{Fano graph of an admissible lattice}
\newcommand{\FNd}{\Fn^d}

All lattices
considered in this note are even.
We
maintain the notation consistent with  \cite{degt:lines}:
\roster*
\item
$\bA_p$, $p\ge1$, $\bD_q$, $q\ge4$, $\bE_6$, $\bE_7$, $\bE_8$
are the \emph{negative definite} root
lattices generated by the indecomposable root systems of the same name
(see~\cite{Bourbaki:Lie});
\item
$\bU:=[0,1,0]$ is 
the unimodular even lattice of rank~$2$;
\item $L(n)$ denotes the lattice obtained by the scaling of
a given lattice~$L$ by a fixed integer  $n\in\Z$
(occasionally, we use $n\in\Q$); this is opposed to
$nL$ that stands for the orthogonal direct sum of  $n$ copies of the lattice $L$;
\item $\bL:=2\bE_8\oplus3\bU$ stands for the \emph{\KK-lattice};
\item
for a lattice $L$, $L\dual:=\Hom(L,\Z)$ denotes the dual lattice;
\item
if $L$ is non-degenerate, $\discr L:=L\dual\!/L$ is the
\emph{discriminant group}, see~\cite{Nikulin:forms};
it inherits from~$L\dual$ a $\Q/2\Z$-valued quadratic form~$q$
(see \autoref{sec.primiso} below);
\item
the  inertia indices of the quadratic form
$L \otimes \R$ are denoted by $\sigma_{\pm,0}(L)$;
\item a
lattice~$S$ with $\Gs_+S=1$
is called \emph{hyperbolic}.
\endroster
Furthermore, we follow the standard lattice--theoretic notation as in,
\eg,~\cite{Conway.Sloane,Nikulin:forms}.

Recall that a \emph{{$\ttD$}-polarized lattice} $S\ni h$ is defined as a
non-degenerate hyperbolic lattice~$S$ equipped with a
distinguished vector $h$ of positive square $h^2 = \ttD$,
called the \emph{degree}.
To simplify our exposition (\cf. \cite[Lemma~2.27.3]{Degtyarev.Rams}),
we assume that
\begin{equation*} 
\ooD > 4 \, .
\end{equation*}
For a polarized lattice $S\ni h$ and $d\in\N$, we define
\[*
\root_d(S,h):=\bigl\{r\in S\bigm|\text{$r^2=-2$, $r\cdot h=d$}\bigr\}.
\]

\definition[\cf. {\cite[\S\,2]{Degtyarev.Rams}}] \label{def-adm-geom}
A polarized lattice $S\ni h$ is called \emph{admissible} if
none of the vectors $e \in S$ satisfies
either of the following conditions
\roster
\item\label{e.h=0}
$e^2 = -2$  and  $e \cdot h = 0$  (\emph{exceptional divisor});
\item\label{e.h=2}
$e^2 = 0$  and  $e \cdot h = 2$  (\emph{$2$-isotropic vector}).
\endroster
An admissible lattice is \emph{geometric}
if it admits
a primitive
isometry to~$\bL$.

More restrictively, given an admissible (resp.\ geometric) lattice $S\ni h$
and a subset $\graph\subset\root_d(S,h)$, we say that $S$ is
\emph{$\graph$-admissible} (resp.\ \emph{$\graph$-geometric}) if,
for each $1\le r<d$,
there is no vector $e\in\root_r(S,h)$ such that
\roster[\lastitem]
\item\label{e.v<0}
$e\cdot v<0$ for some $v\in\graph$ (\emph{$\graph$-inadmissible root}).
\endroster
\enddefinition

For an admissible lattice $S\ni h$ (so that, in particular,
$\root_0(S,h)=\varnothing$),
we put  $\basis\fp_1 :=  \root_1(S,h)$ and,  for $d > 1$, we define recursively
(\cf. \cite[\S\,2.2]{Degtyarev.Rams})
\[*
\basis\fp_d:=\bigl\{r\in\root_d(S,h)\bigm|
\text{$r\cdot e\ge0$ for all $e\in\basis\fp_r$, $1\le r<d$}\bigr\}.
\]

\definition \label{def-fano}
We  define
the (plain) \emph{degree~$d$ Fano graph} of an admissible lattice
$S\ni h$
as the set of vertices
\[ \label{eq-def-Fano}
\FNd(S,h):= \basis\fp_d,
\]
with two vertices $l_1\ne l_2$ connected by an edge of multiplicity
$l_1\cdot l_2$.
\enddefinition

For example,
the degree~$1$ Fano graph $\Fn^1(S,h)$
is $\Fn_{\emptyset}(S,h){}=\root_1(S,h)$ in the sense of
\cite[(2.6)]{Degtyarev.Rams}.

\newcommand{\Fnr}{\Fn^r}

\remark\label{rem.reflections}
It
follows immediately from the general theory of groups generated by
reflections that any root $e\in\root_d(S,h)$, $d>0$, is a \emph{positive} linear
combination
\[*
e=\sum_in_ie_i,\quad n_i\in\N,\quad
e_i\in\bigcup_{r=0}^{d}\Fnr(S,h).
\]
Geometrically (\cf. the discussion below), this means that any $(-2)$-curve
of positive degree splits into a union of smooth rational components
(clearly, of smaller degrees unless the original curve is smooth).
It follows that
\[*
\Fn_d(S,h)=\bigl\{r\in\root_d(S,h)\bigm|
 \text{$r\cdot e\ge0$ for all $e\in\root_r(S,h)$, $1\le r<d$}\bigr\},
\]
\ie, there is no need to compute the intermediate sets $\Fn_r(S,h)$.
Furthermore, the $\graph$-admissibility condition
in \autoref{def-adm-geom}
merely means that
$\graph\subset\Fn_d(S,h)$.
\endremark

\subsection{Primitive isometries into $\bL$\noaux{ (see~\cite{Nikulin:forms})}} \label{sec.primiso}
Consider a non-degenerate even lattice $N$ and extend the bilinear form to a
$\Q$-valued form on $N\otimes\Q$.
The dual group $N\dual$ can be identified with
the subgroup
\[*
N\dual=\bigl\{x\in N\otimes\Q\bigm|
\text{$x\cdot y\in\Z$ for all $y\in N$}\bigr\}\supset N,
\]
and the \emph{discriminant group}
(denoted~$q_N$ in~\cite{Nikulin:forms}) is defined as the finite quotient
\[*
\discr N:=N\dual\!/N.
\]
This group is equipped with a $(\Q/\Z)$-valued symmetric bilinear
form and its $(\Q/2\Z)$-valued quadratic extension:
\[*
(x\bmod N)\otimes(y\bmod N)\mapsto x\cdot y\bmod\Z,\quad
x\bmod N\mapsto x^2\bmod2\Z.
\]
We often abbreviate
\[*
\CN:=\discr N\quad\text{and}\quad
\CN_p:=\discr_pN:=\CN\otimes\Z_p\ \text{for a prime~$p$}
\]
and use the shorthand
\[*
\cl[\xi]:=\xi\bmod N\in\CN
\]
for the discriminant class of a vector $\xi\in N\dual$.
We also use the notation $\<x\>\subset\CN$ for the subgroup generated by
$x\in\CN$, and we abbreviate $\<\xi\>:=\<\cl[\xi]\>$ for $\xi\in N\dual$.

Obviously, $\CN=\bigoplus_p\CN_p$ over all primes~$p$.
The discriminant form on~$\CN$ is non-degenerate in
the sense that the map
\[*
\CN\to\Hom(\CN,\Q/\Z),\quad
x\mapsto(y\mapsto x\cdot y\bmod\Z),
\]
is an isomorphism.
Besides, one has $\ls|\discr N|=\ls|\det N|$.

\remark \label{rem:compdiscr}
To find the group $\discr_p N$, we
indicate a subset
therein,
whereupon it suffices to
compute the Gram matrix~$G$ and check that
$\ppart{p}(\det G\cdot\det N)=0$
(\cf. \autoref{s.12A1} below):
the non-degeneracy would imply that the set is a basis.
\endremark


For $p=2$, we distinguish between \emph{even} (\ie, such that $x^2=0\bmod\Z$
for each element
$x\in\CN_2$ \emph{of order~$2$}) and \emph{odd} (otherwise) forms. In other
words, the form $\CN_2$ is odd if and only if it has an element of square
$\pm\frac12\bmod2\Z$ (necessarily of order~$2$ and necessarily an orthogonal
direct summand).

The \emph{length} $\ell(\CN)$ of a finite abelian group~$\CN$ is the minimal
number of generators. We usually describe a discriminant form by its
\emph{Gram matrix} $[a_{ij}]$ in a minimal basis
$x_1,x_2,\ldots$\,:
the
diagonal entries $a_{ii}:=x_{ii}^2$ are defined modulo~$2\Z$, whereas all
other entries $a_{ij}:=x_i\cdot x_j$, $i\ne j$, are only defined modulo~$\Z$.

If $p>2$ or $p=2$ and $\CN_2$ is even, one can speak about the
\emph{determinant}
$\det\CN_p$:
\[*
\ls|\CN_p|\det\CN_p\in\Z_p\units,\quad
\text{well defined modulo $(\Z_p\units)^2$}.
\]
In~\cite{Nikulin:forms}, this is defined as the determinant (divided by $\ls|\CN_p|^2$,
which is irrelevant modulo squares) of the only $\Z_p$-lattice~$K_p$ such that
\[*
\rank K_p=\ell(\CN_p),\quad \discr K_p\simeq\CN_p.
\]
(Such a lattice
exists and,
unless $p=2$ and $\CN_2$ is odd,
is unique up to isomorphism.)
Alternatively, this is just the determinant of the
Gram matrix in a carefully chosen minimal basis for $\CN_2$: \eg, it
suffices to require that any two basis elements of distinct orders must be
orthogonal. Indeed, such a basis immediately gives us~$K_p$: we merely change
all $p^{-k}$ to $p^k$, $k\in\N$, in the matrix.

If $\CN_2$ is odd, then $\det\CN_2$ is only defined modulo $5\in\Z_2\units$,
making it redundant in \autoref{th.Nikulin} below.

The determinant is multiplicative,
$\det(\CN\oplus\Cal M)=(\det\CN)(\det\Cal M)$.
Each non-degenerate finite
quadratic form splits into orthogonal direct sum of the following elementary
blocks: for a prime $p>2$,
\[*
\Cal{W}=\biggl[\frac{2m}{p^k}\biggr],\quad
\text{$\gcd(m,p)=1$ and only
	$p^k\det\Cal{W}=\biggl(\frac{2m}{p}\biggr)\in\{\pm1\}$ matters},
\]
and for $p=2$,
\[*
\Cal{W}=\biggl[\frac{m}{2^k}\biggr],\quad
\text{$m$ is odd and only
	$2^k\det\Cal{W}=m\bmod8\in\{\pm1,\pm5\}$ matters},
\]
or
\[*
\CU_{2^k}:=\bmatrix0&1/2^k\\1/2^k&0\endbmatrix,\qquad
\CV_{2^k}:=\bmatrix1/2^{k-1}&1/2^k\\1/2^k&1/2^{k-1}\endbmatrix;
\]
we have
\[*
2^{2k}\det\CN=
\begin{cases}
	-1\bmod8,&\text{if $\CN=\CU_{2^k}$},\\
	\phantom-3\bmod8,&\text{if $\CN=\CV_{2^k}$}.
\end{cases}.
\]

\theorem\label{th.Nikulin}
{\rm (see \cite[Theorem~1.12.2]{Nikulin:forms})} A non-degenerate hyperbolic lattice $N$ admits a primitive isometry into
the \emph{$K3$-lattice} $\bL:=2\bE_8\oplus3\bU$ if and only if
the following conditions are satisfied simultaneously\rom:
\roster
\item\label{Nikulin.r}
one has $r:=22-\rank N\ge2$\rom;
\item\label{Nikulin.p}
for each prime $p>2$, either
\roster*
\item
$\ell(\CN_p)\le r-1$, or
\item
$\ell(\CN_p)=r$ and $\ls|\CN|\det\CN_p=1\bmod(\Z_p\units)^2$\rom;
\endroster
\item\label{Nikulin.2}
for $p=2$, either
\roster*
\item
$\ell(\CN_2)\le r-2$, or
\item
$\ell(\CN_2)=r$ and the form $\CN_2$ is odd, or
\item
$\ell(\CN_2)=r$, the form $\CN_2$ is even, and
$\ls|\CN|\det\CN_2=\pm1\bmod8$.
\done
\endroster
\endroster
\endtheorem


\subsection{Lattice associated to a graph} \label{sec.lattice-of-graph}
Let $\graph$ be a loop-free (multi-)graph.
Then, as in  \cite{degt:lines}, we put  $\Z\graph$ to denote
the lattice freely generated by the
vertices of~$\graph$,
with the bilinear form defined by
$v^2=-2$ and $u\cdot v=r$
for $u,v\in\graph$
connected by an $r$-fold edge (no edge meaning a $0$-fold edge).

Recall that $\graph$ is called
\emph{elliptic} (resp. \emph{parabolic}, resp. \emph{hyperbolic})
iff  $\Z\graph$ is negative definite (resp.  $\Gs_+(\ZZ \graph) = 0$ and  $\sigma_0(\ZZ \graph) >0$, resp.  $\Gs_+(\ZZ \graph) = 1$).

Given a pair of integers $\ooD,d>0$, we mimic
\cite[\S\,2.4]{degt:lines} and associate to the graph $\graph$  the $2\ooD$-polarized lattice
\[ \label{eq:fanolattice}
\Fano_{\ttD}^d(\graph):=(\Z\graph+\Z h)/\ker,\qquad h^2=\ttD,\quad
\text{$h\cdot v=d$ for $v\in\graph$},
\]
where $\ker:=\ker(\Z\graph+\Z h)=(\Z\graph+\Z h)^\perp$
is the
radical
of the bilinear form.

\remark\label{rem.colored}
More
generally, one can consider a colored graph with the
coloring function $\deg\:\graph\to\{0,\ldots,d\}$ and let $h\cdot v=\deg(v)$
in~\eqref{eq:fanolattice}.
Most definitions below extend to this case in an obvious way.
In this paper, most of the time we assume that
$\deg(v)=d$ is a constant
fixed in advance.
\endremark

By \cite[Proposition 1.4.1]{Nikulin:forms},
there is a
bijection
between
the (isomorphism classes of)
even finite index extensions
$N'\supset N$
and isotropic subgroups (\latin{a.k.a}.\ \emph{kernels})
$\CK\subset\CN$: one has
\[*
\CK:=N'/N\subset\discr N,\qquad \discr N'=\CK^\perp\!/\CK.
\]
Thus, one can pick a
kernel
$\CK\subset\discr \Fano_{\ttD}^d(\graph)$ and consider the extension
$\Fano_{\ttD}^d(\graph,\CK)$ of $\Fano_{\ttD}^d(\graph)$ by~$\CK$.
The kernel~$\CK$ is called \emph{$\graph$-admissible} if so is
$\Fano_{\ttD}^d(\graph,\CK)$.

\definition \label{def-adm-graph}
We call a  graph $\graph$ \emph{admissible}
in degree~$\ttD$ (with respect to $d$) if the lattice
$\Fano_{\ttD}^d(\graph)$ is $\graph$-admissible.
An admissible graph~$\graph$ is called \emph{subgeometric}
if there is a $\graph$-geometric finite index extension
$\Fano_{\ttD}^d(\graph,\CK)$;
the kernel $\CK$ of such an extension
is called a \emph{geometric kernel} for~$\graph$.
If such a geometric kernel~$\CK$ can be chosen so that
$\graph=\Fn_d\Fano_\ttD^d(\graph,\CK)$, then the graph~$\graph$ is called
\emph{geometric}.

Occasionally, we speak about a \emph{partial} or, more precisely,
\emph{$(\bmod\,p)$-geometric kernel}
$\CK_p\subset\discr_p\Fano_{\ttD}^d(\graph)$, meaning that
$\CK_p$ is $\graph$-admissible and
the $p$-torsion part $(\CK_p^\perp/\CK_p)_p$ satisfies
respective condition~\iref{Nikulin.p}
or~\iref{Nikulin.2} in \autoref{th.Nikulin}.
\enddefinition



Let
$\Sigma$ be a connected parabolic graph (\latin{a.k.a}.\ affine Dynkin
diagram). The radical $\ker\Z\Sigma$ has a unique minimal positive generator
\[
k_\Sigma=\sum_{v\in\Sigma}n_vv,\quad n_v\in\Z_+,
\label{eq.k}
\]
and the \emph{degree} $\deg(\Sigma)$ is defined as the coefficient sum
$\sum n_v$. We use the notation
$\p\Sigma\subset\Sigma$ for an associated (elliptic) Dynkin diagram
obtained from~$\Sigma$ by removing (any) one vertex $v\in\Sigma$ with
$n_v=1$ in \eqref{eq.k}.
(This operation is not quite a functor as a vertex~$v$ that can be removed is
typically not unique.) In general, for a parabolic graph~$\graph$,
we denote by $\p\graph$ the graph obtained by
replacing each parabolic component~$\Sigma$ with $\p\Sigma$
(and keeping all elliptic components intact).

If $\graph$ is a parabolic graph, one has
$\Gs_+(\Z\graph+\Z h)=1$ if and only if all parabolic components of $\graph$
have the same degree; this common value $\deg(\graph)$ is called the
\emph{degree} of~$\graph$. In this case, we can define
\[
k:=k_\Sigma\in\Fano_\ttD^d(\graph),\qquad
\kk:=\frac{1}{\deg(\Sigma)}k_\Sigma\in \Fano_{\ttD}^d(\graph)\dual,
\label{eq:defkkappa}
\]
where $\Sigma\subset\graph$ is any parabolic component:
all images in $\Fano_{\ttD}^d(\graph)\dual$ are equal.


\lemma\label{lem.no.k}
Let $\graph$ be a parabolic graph. Then any
kernel $\<qk\>\subset\discr\Fano_\ttD^d(\graph)$, $q\in\Q\sminus\Z$,
contains an exceptional divisor or a $\graph$-inadmissible root.
\endlemma

\proof
Pick a parabolic component $\Sigma\subset\graph$ and let
$f:=\deg(\Sigma)=\deg(\graph)$. We recall the following less known property
of the elliptic Dynkin diagram~$\p\Sigma$: the set of values taken by the
coefficient sum $\sum n_e$ over all positive roots $r=\sum n_vv$,
$v\in\p\Sigma$, $n_v\in\N$, is precisely the set $\{1,\ldots,f-1\}$.

We can assume that $0<q\le1/2$
(and, clearly, $fdq\in\Z$). Then $0<fq\le f-1$ and, letting
$m:=\lceil fq\rceil$, we can find a positive root $r\in\Z\p\Sigma$
with the coefficient sum~$m$. As another property of Dynkin diagrams,
$r\cdot v<0$ for some $v\in\p\Sigma$. Thus, $\cl[-qk]$ contains the root
$e:=-qk+r$ of degree $0\le e\cdot h=(m-qf)d<d$ and such that $e\cdot v<0$ for
some $v\in\graph$.
\endproof

The following statement would explain the
$d\deg(\graph)$-periodicity of the degree~$d$ curve counts, but it is not strong
enough; therefore, we omit the proof.

\proposition[\cf. {\cite[Proposition~3.7]{degt:lines}}] \label{prop:periodic}
Let $\graph$ be a parabolic graph, and let $f := \deg(\graph)$. If
$\ttD',\ttD'' \geq (df)^2 + df$
and $\ttD' = \ttD'' \bmod df$, then $\graph$ is geometric in degree
$\ttD'$
if and only if so it is
in degree $\ttD''$.
\endproposition


\lemma\label{lem.k.perp}
Let
$\graph$ be a parabolic graph and $d>0$ a fixed integer. Then, for any
integer $\ooD\gg d$, any vector $e\in\Fano_{\ttD}^d(\graph)\dual$ that is
either
\roster
\item\label{k.perp.exceptional}
an exceptional divisor, as in \autoref{def-adm-geom}\iref{e.h=0}, or
\item\label{k.perp.isotropic}
a $2$-isotropic vector, as in \autoref{def-adm-geom}\iref{e.h=2}, or
\item\label{k.perp.root}
more generally, an \emph{$r$-root}, \ie, $e^2=-2$ and $r:=e\cdot h\le d$
\endroster
lies in $k^\perp$. In case~\iref{k.perp.isotropic}, this implies that
$e=2\kk/d$ is
a multiple of~$k$.

Conversely,
any root $e\in k^\perp\subset\Fano_{\ttD}^d(\graph)\dual$ such
that $0\le r:=e\cdot h<d$ generates a $\graph$-inadmissible finite index
extension of $\Fano_{\ttD}^d(\graph)$.
\endlemma

{Essentially, in items~\iref{k.perp.exceptional} and~\iref{k.perp.root} we
assert that any smooth rational curve of low degree is a fiber component of
the elliptic pencil in question, which is well known.
Though, the presence of such a curve
most likely means
that the pencil is not what we think it is, \cf. \autoref{ex.low.degree}
below.
}

\proof
We compute the Gram matrix $G$ of $h,k,e$ and assert that $\det G\ge0$, since
the lattice is hyperbolic. In cases~\iref{k.perp.exceptional} (letting $r=0$)
and~\iref{k.perp.root} this results in
\[*
2x(dfr - \ooD x) + 2(df)^2\ge0,
\]
where $x:=e\cdot k\in\Z$;
in case~\iref{k.perp.isotropic}, we obtain
\[*
2x(2df-\ooD x)\ge0.
\]
In all cases, since $x\in\Z$, we conclude that $x=0$ for $\ooD\gg0$.

For the converse statement,
consider the root lattice
\[*
\Z\p\graph=\Z\graph/\ker;
\]
in this lattice, the vertices of~$\p\graph$ constitute a
standard Dynkin basis (\latin{a.k.a}.\ Weyl chamber).
By \autoref{lem.no.k}, we can assume that $e\notin\cl[qk]$, $q\in\Q$;
hence, the root $\bar{e}:=e\bmod(\ker\Z\graph)\in(\Z\p\graph)\dual$ represents a
non-trivial discriminant class.
The root lattice
\[*
(\Z\p\graph+\Z\bar{e})/\ker\supset\Z\p\graph
\]
is a proper finite index extension;
hence, it has strictly smaller Weyl chambers, {\4 so} it cannot be
$\graph$-admissible.
\endproof


\remark \label{rem:smallperiod}
Let
$\graph$ be a parabolic graph and $\ooD\gg0$.
\autoref{lem.k.perp} implies that any exceptional divisor ($r=0$) or
degree~$r$
root
$e\in\Fano_{\ttD}^d(\graph)\dual$ is of the form
\[*
e=\sum_{v\in\graph} q_vv,\quad \text{where}\
q_v\in\Q,\quad \sum_{v\in\graph}q_v=\frac{r}{d},\quad
e^2=-2,\quad
e\cdot u\in\Z\ \text{for all $u\in\graph$}.
\]
\latin{De facto}, we are working with the discriminant groups of root
systems, which are very well understood. In particular,
we can speak about
shortest representatives of the discriminant classes.
It is this fact and the periodicity of
the relevant part of the discriminant (rather than that of the lattice
itself)
that explains
period $\deg(\Sigma)$ (rather than $d\deg(\Sigma)$ given by
\autoref{prop:periodic}) in most statements.
\endremark

\definition\label{def.ex.orbit}
With \autoref{lem.k.perp} in mind, an isotropic class
$\Ga\in\discr\Fano_\ttD^d(\graph)$ is called \emph{exceptional} if it is
represented by a root $e\in\Fano_\ttD^d(\graph)\dual$ orthogonal
to~$k$ and such that $0\le e\cdot h<d$.
According to \autoref{rem:smallperiod},
such classes are independent of~$\ooD$.
An isotropic subgroup $\CK\subset\discr\Fano_\ttD^d(\graph)$ is
called
\emph{exceptional} if it contains an exceptional class. This
terminology extends, in the obvious way, to the $(\Aut\graph)$-orbits of
isotropic classes/subgroups. We can also speak about \emph{minimal}
exceptional subgroups and \emph{minimal} exceptional
classes as those generating such subgroups.
\enddefinition


\subsection{Degree~$d$ smooth rational curves on \smooth $K3$-surfaces}\label{s.preli.geom}
A \emph{degree~$\ttD$} polarized $K3$-surface is a pair $(X,h)$, where $X$
is a $K3$-surface and $h\in\NS(X)$ is a very ample divisor, $h^2=\ttD$;
\cf. \cite{degt:lines,rams.schuett:24,rams.schuett.smoothodd}.
Usually, the polarization~$h$ is considered part of the structure and omitted
from the notation.

We imitate \cite{degt:lines,rams.schuett:24}  (\cf. also \cite{Degtyarev.Rams,rams.schuett:quasi})
and define the Fano graph of smooth rational degree~$d$ curves as the
set
\[
\FNd(X,\ttH)
:= \{\text{smooth rational curves $C\subset X$ with $C\cdot\ttH = d$}\},
\]
where two vertices $C,C'$ are connected by an edge of multiplicity
$C\cdot C'$.
We focus on curves of a fixed degree $d$;
hence, contrary to
\cite{Degtyarev.Rams,rams.schuett:24}, we treat
$\FNd(X,\ttH)$ as a plain graph;
in other words, we assign a fixed color~$d$ to all vertices.

We say that the N\'{e}ron--Severi lattice $\NS(X)$ is
\emph{spanned by degree~$d$ curves}
if it is
a finite index extension
of its sublattice generated by the classes of
\emph{smooth rational} degree~$d$ curves  on~$X$
\emph{and the polarization~$\ttH$}.

Torelli theorem for \KK-surfaces and the theory of lattice-polarized \KK-surfaces immediately yield the following version of
\cite[Theorem~3.9]{Degtyarev.Rams}

\theorem \label{thm:K3}
Let $d>0$ and
$\ooD>4$ be  fixed integers. A graph $\graph$ is geometric in
degree~$\ttD$ with respect to $d$ if and only if
$\graph\sameg\FNd(X,\ttH)$ for a
degree~$\ttD$ polarized \KK-surface $(X,\ttH)$
such that $\NS(X)$ is \qgen by
degree~$d$ curves.
\done
\endtheorem


\remark\label{rem.Schur}
If $d>1$,
the realizability criterion for a graph of degree~$d$ curves
\emph{must} appeal to curves of lower degree. For example, famous Schur's
quartic~\cite{Schur:quartics} has (obviously admissible) N\'{e}ron--Severi lattice
spanned by conics, but most of these conics are reducible,
see~\cite{degtyarev24conicsquartics}.
In the parabolic context, see \autoref{ex.low.degree} below.

In particular, to lift the condition that $\NS(X)$ should be spanned by
$d$-curves, one would have to consider the full Fano graph of smooth rational
curves of degree up to~$d$, \cf. \autoref{rem.colored} and
\eqref{eq.not-spanned} below.
Clearly, the requirement that the graph should be
\emph{sub}-geometric is necessary for its realizability in any case.

\endremark

\autoref{thm:K3} justifies the fact
that, once
the integers $d$ and $\ooD$ have been fixed, one can
identify geometric graphs
with configurations of smooth rational degree~$d$ curves on \KK-surfaces. In
particular, when we speak of a \emph{pencil} (resp.\ its \emph{fiber}),
we mean a parabolic
graph (resp.\ connected parabolic subgraph).


\subsection{Simplest applications: type $\ttA1$ pencils}\label{S.A1}

We
start with a simple example illustrating the machinery described in this
section.

 \example \label{s.12A1}
Fix an \emph{odd} integer $d>1$.
We claim that,
for an \emph{odd} integer $\ooD \gg d$,
 there is no \KK-surface of degree $\ttD$
with
 $\FNd(X,\ttH) \simeq  12\tA_1$.
 By  \autoref{thm:K3} it suffices to prove that
 \begin{equation}
 	\graph := 12\tA_1 \mbox{ fails to be subgeometric in degree $\ttD$
 		with respect to  $d$.}	
 \end{equation}
A basis for $\Fano_{\ttD}^d(\graph)$ is
 \[*
 h,k,a_1,\ldots,a_{12},
 \]
 where
$a_i,a_i'$, $i=1,\ldots,12$, are the fiber components and $k=a_i+a_i'$ is as
in~\eqref{eq.k}.
Obviously,
 \[*
 \det \Fano_{\ttD}^d(\graph)  =-2^{14}d^2,
 \]
as can be seen from the orthogonal decomposition over $\Q$.
Then, by \autoref{rem:compdiscr},
$\discr_2\Fano_{\ttD}^d(\graph)\simeq(\Z/2)^{10}\oplus(\Z/4)^2$
(recall that $d$ is odd)
is generated by
 \[
 \Ga_i:=\frac12(a_i-a_{i+1}),\ i=1,\ldots,10,\quad
 \eta:=\frac12h+\frac14\sum_{i=1}^{12}a_i,\quad
 \kappa_4:=\frac14k-\frac12a_{12}.
 \label{eq.basis.12I1}
 \]
with the Gram matrix
 \[
 \bA_{10}\bigl(\tfrac12\bigr)\oplus
 \bmatrix s&\phantom-1/4\\1/4&-1/2\endbmatrix,\quad
 s :=3d+\frac12\ooD-\frac32.
 \label{eq.Gram.12I1}
 \]
Thus,
 \[*
 \ell(\discr_2  \Fano_{\ttD}^d(\graph))+\rank  \Fano_{\ttD}^d(\graph)=26=\rank(\bL)+4,
 \]
\ie,
we need a finite index extension of $\Fano_{\ttD}^d(\graph)$ by at least two
isotropic vectors (see \autoref{th.Nikulin}).
The $\OG_h\bigl(\Fano_{\ttD}^d(\graph)\bigr)$-orbits of isotropic vectors in $\discr_2 \Fano_{\ttD}^d(\graph)$ are
 represented by
 \[*
 \gathered
 \cl[2\kappa_4]=\cl[\kk]\ni\frac12(a_{12}'-a_{12}),\qquad
 \cl[\Ga_1+\Ga_3]\ni\frac12(a_1-a_2+a_3-a_4),\\
 \cl[\Ga_1+\Ga_3+\Ga_5+\Ga_7]\ni\frac12\sum_{i=1}^{8}a_i,\qquad
 \cl[2\eta]\ni\frac12\sum_{i=1}^{12}a_i;
 \endgathered
 \]
Since $\ttD:= h^2=4m+2$, we have $s=3d+m-1\in\Z$ in~\eqref{eq.Gram.12I1};
hence, there also is
an orbit of order~$4$
classes
$\cl[\lambda]$ such that $\cl[2\lambda]=\cl[2\eta]$.
The first two orbits are
exceptional.
Now, it is straightforward that, up to $\OG_h\bigl(\Fano_{\ttD}^d(\graph)\bigr)$,
 the only admissible index~$4$ extension $N'\supset \Fano_{\ttD}^d(\graph)$ that \emph{might} embed
 to the $K3$-lattice $\bL$ is that by
 \[
 \frac12\sum_{i=1}^{8}a_i\quad\text{and}\quad \frac12\sum_{i=5}^{12}a_i,
 \label{eq.kernel.12I1}
 \]
 with the new discriminant
 \[
 \discr_2 N'=\CK^\perp\!/\CK\sameg\CV_2^3\oplus\begin{cases}
 	\CV_4, & \text{if $h^2=0\bmod4$}, \\
 	\CU_4, & \mbox{if $h^2=2\bmod4$}
 \end{cases}
 \label{eq.new-discr}
 \]
 of length~$8=\rank\bL-\rank\Fano_{\ttD}^d(\graph)$,
\cf. \autoref{rem.new-discr} below.
The last case $\CV_2^3\oplus\CU_4$
with $\ooD$ odd
violates Nikulin's existence criterion,
\viz.
\autoref{th.Nikulin}\iref{Nikulin.2}.
\qed
\endexample

\remark\label{rem.new-discr}
The
isomorphism type of the form~\eqref{eq.new-discr}
the group
 $(\Z/2)^6\oplus(\Z/4)^2$ is uniquely determined by its Brown invariant
 $\Br=4\bmod8$ and the number of order~$4$ classes of square
 $0\bmod\Z$; clearly, the latter gets divided by~$4$ when passing to
 $\CK^\perp\!/\CK$.
However, in the sequel we choose to refrain from the discussion of various shortcuts
that can be used: all forms can be computed explicitly.
\endremark



Thanks to~\cite{rams.schuett.smoothodd}, we are not concerned with the
\emph{existence} of type~$\tA_1$ pencil for $d$ or~$\ooD$ even. Therefore, we
confine ourselves to the following negative statement.

\proposition \label{lem:A1}
Let
$h\gg d>1$
be \emph{odd} integers.
Then
neither
of the
graphs
$$
12\tA_1,  \, 11\tA_1
$$
is subgeometric in degree $\ttD$ with respect to $d$.
\endproposition


\proof
It suffices to
imitate  the  argument
of \autoref{s.12A1}
to show that
the graph $\graph:= 11\tA_1$
fails to be subgeometric.
This time,
$\discr_2\Fano_{\ttD}^d(\graph)\sameg(\Z/2)^{10}\oplus(\Z/8)$
is generated by $\Ga_1,\ldots,\Ga_{10}$ as
in~\eqref{eq.basis.12I1} and
\[*
\bar\eta:=\frac12h+\frac18k+\frac14\sum_{i=1}^{11}a_i,\quad
\bar\eta\cdot\Ga_i=0,\quad
\bar\eta^2=3d+\frac12\ooD-\frac{11}{8}.
\]
Up to
$\OG_h\bigl(\Fano_{\ttD}^d(\graph)\bigr)$,
the only admissible extension $N'\supset\Fano_{\ttD}^d(\graph)$ is that of
index~$2$ by the first vector in~\eqref{eq.kernel.12I1}. We have
\[*
\discr_2 N'=\CU_2^3\oplus\CV_2\oplus
\bigl[\bar\eta^2\bigl],
\]
which fails to satisfy
\autoref{th.Nikulin}\iref{Nikulin.2}
unless $2\ooD = 0\bmod4$.
\endproof


\section{Type $\ttA3$ pencils}\label{S.A3}

In this section, we also assume that \emph{$d$ is odd} and mainly prove the
\emph{non}-existence of \KK-surfaces of degree $2\ooD=2\bmod4$
with several $\tA_3$-configurations
of smooth rational curves of odd degree $d$.
We index the $\tA_3$ type components of the graph in question and denote by
$a^0_i,a^1_i,a^2_i,a^3_i$ the vertices
(numbered consecutively along the cycle $\tA_3$)
of the $i$-th component.
We also use the polarization~$h$ and
isotropic vector $k=a^0_i+\ldots+a^3_i=\const(i)$,
see~\eqref{eq.k}.

\subsection{The pencil $6\ttA3$}\label{s.6A3}
The following proposition holds.

\proposition \label{lemma.6tA3}
Let
$\ooD\gg d>1$
be \emph{odd} integers.
Then the
graph  $\graph:=6\tA_3$
fails to be  subgeometric in degree $\ttD$ with respect to $d$.
\endproposition
\proof
For the basis of $\Fano_{\ttD}^d(\graph)$,  we  take
\[*
h,k,a_i^r,\ i=1,\ldots,6,\ r=1,2,3.
\]
Obviously, we have $\det\Fano_{\ttD}^d(\graph)=-2^{16}d^2$ and, letting
\[*
\ba_i:=a_i^1+2a_i^2+3a_i^3,\quad
\sa_i:=3a_i^1+4a_i^2+3a_i^3,\quad
i=1,\ldots,6,
\]
the group $\discr_2\Fano_{\ttD}^d(\graph)   \sameg(Z/2)^2\oplus(\Z/4)^4\oplus(\Z/8)^2$ is generated
by
\[
\gathered
\eta_2:=\frac12h+\frac12\sum_{i=1}^{6}(a_i^1+a_i^2),\quad
\Gb_6:=\frac12(a_6^1-a_6^3),\\
\Ga_i:=\frac14(\ba_i-\ba_{i+1}),\ i=1,\ldots,4,\\
\eta_8:=\frac{d}{4}h+
\frac{1-\ooD}{8}k
+\frac18\sum_{i=1}^{6}\sa_i,\quad
\kappa_8:=\frac18k+\frac14\sum_{i=1}^{5}\ba_i.
\endgathered
\label{eq.6A3.basis}
\]
In this basis, the matrix of the quadratic form
is
\[
\bmatrix s_2&1/2\\1/2&1\endbmatrix\oplus\bA_4\bigl(\tfrac14\bigr)\oplus
\bmatrix s_8&1/8\\1/8&1/4\endbmatrix,\quad
s_2:=\frac12\ooD - 1,\quad
s_8:=\frac18-\frac{1}{8}\ooD d^2.
\label{eq.6A3.form}
\]
Thus, we find that the form $\discr_2\Fano_{\ttD}^d(\graph)$
depends on the residue
$\ooD\bmod4$
only.
(For
$n$ odd, the apparent dependence of $s_8$
on the
\emph{odd} residue $d\bmod8$
is irrelevant due to the obvious isomorphism of the last summand  to
$\CU_8$.)
By \autoref{lem.k.perp},
neither~$\ooD$
nor~$d$ affect the minimal exceptional classes;
up to
$\OG_h\bigl(\Fano_{\ttD}^d(\graph)\bigr)$
they are
\[
\cl[4\kappa_8]=\cl[2\kk]\ni\frac12(a_1^0+a_1^1-a_1^2-a_1^3),\qquad
\cl[2\Ga_1]\ni\frac12(a_1^1+a_1^3-a_2^1-a_2^3).
\label{eq.6A3.classes}
\]
We conclude that, for
$\ooD\gg d$,
the existence of {\4 the} $K3$-surface
does not depend on~$d$
provided that the latter is odd; hence, as in the case $d=1$ of lines
(see~\cite{degt:lines}),
$\graph$
cannot be subgeometric
unless  $\ttD=4\bmod8$.
\endproof

\remark \label{rem:existence6tA3}
{We do not discuss
the \emph{existence} of surfaces with a pencil $6\tA_3$
of degree~$d$ curves as we do not need it for the proof of
\autoref{thm.main}. It appears that, for $d$ odd, apart from $\ttD=4\bmod8$
(\cf. \cite[\S\,11.1]{rams.schuett:24}) one should also impose a
$q$-balancing condition,
\cf. \autoref{lem.q} and \autoref{th.4A4+2A1} below.
}
\endremark

\subsection{The pencil $5\ttA3+2\bA\sb1$}\label{s.5A3+2A1}
Below we prove the following statement.

\proposition \label{lemma.5tA3plus2bA1}
Let
$\ooD\gg d>1$
be \emph{odd} integers.
Then the
graph $\graph:=5\tA_3+2\bA_1$
fails to be  subgeometric in degree $\ttD$ with respect to $d$.
\endproposition

\proof
We consider the following basis of $\Fano_{\ttD}^d(\graph)$:
\[*
h,k,a_i^r,b_1,b_2,\ i=1,\ldots,6,\ r=1,2,3,
\]
where $b_1,b_2$ are the two $\bA_1$ fibers.
The 
group
$\discr_2 \Fano_{\ttD}^d(\graph)  \same(\Z/4)^5\oplus(\Z/8)^2$ is generated by
$\Ga_1,\ldots,\Ga_4$ as in~\eqref{eq.6A3.basis} and
\[*
\gathered
\bar\kappa_4:=\frac14k+\frac14\sum_{i=1}^{5}\ba_i+\frac12b_1,\quad
\bar\kappa_8:=\frac78k+\frac12\sum_{i=1}^{5}\ba_i+\frac32b_1,\\
\bar\eta_8:=\frac{d}{4}h+
\frac{3-\ooD}{8}k
 +\frac18\sum_{i=1}^{5}(5a_i^1+a_i^3)+\frac58b_1+\frac18b_2;
\endgathered
\]
the reduced intersection matrix is
\[
\bA_4\bigl(\tfrac14\bigr)\oplus\biggl[\frac74\biggr]\oplus
 \bmatrix1/2&3/4\\3/4&s_8\endbmatrix,
\label{eq.5A3+A2.form}
\]
where $s_8$ is as in~\eqref{eq.6A3.form}. As in \autoref{s.6A3},
$\discr_2\Fano_{\ttD}^d(\graph)$ depends on
$\ooD\bmod4$
only.
The minimal exceptional classes are $\cl[4\bar\kappa_8]=\cl[2\kk]$ and
$\cl[2\Ga_1]$ as in \eqref{eq.6A3.classes} and
\[*
\cl[4\bar\eta_8-2\Ga_1-2\Ga_3]\ni\frac12(a_5^r+a_5^{r+2}-b_1-b_2),
\]
where $r=0$ or $1$ for
$\ooD$ even or odd,
respectively.
Thus,
as in the case $d=1$ of lines, see~\cite{degt:lines},
$\graph$ is not subgeometric unless $\ooD$ is even.
\endproof


\subsection{The pencil $5\ttA3+\bA\sb2$}\label{s.5A3+A2}
The last  $\tA_3$ pencil to
be considered
is  $5\tA_3+\bA_2$.

\proposition \label{lemma.5tA3plusbA2}
Let
$\ooD\gg d>1$
be \emph{odd} integers.
Then the
graph $\graph:=5\tA_3+\bA_2$
fails to be  subgeometric in degree $\ttD$ with respect to $d$.
\endproposition

\proof
We take for the basis
\[ \label{eq-5A3basis}
h,k,a_i^r,b_1,b_2,\ i=1,\ldots,6,\ r=1,2,3,
\]
where $b_1,b_2$ are the two curves in the $\bA_2$ fiber.
The discriminant group
$$
\discr_2\Fano_{\ttD}^d(\graph)\sameg(\Z/2)^2\oplus(\Z/4)^4\oplus(\Z/16)
$$
is generated by
$\Ga_1,\ldots,\Ga_4$ as in~\eqref{eq.6A3.basis} and
\[*
\gathered
\bar\eta_2:=\frac12h+\frac12\sum_{i=1}^{5}(a_i^1+a_i^2)+\frac12(b_1+b_2),\quad
\Gb_1:=\frac12(a_1^1-a_1^3),\\
\bar\eta_{16}:=\frac{d}{4}h-
 \frac{1+\ttD}{16}k
 +\frac18\sum_{i=1}^{5}(5a_i^1+a_i^2)+\frac{5}{4}(b_1+b_2).
\endgathered
\]
The reduced intersection matrix in the basis
$\bar\eta_2,\Gb_1,\Ga_1,\ldots,\Ga_4,\bar\eta_{16}$ is
\[*
\bmatrix
s_2&  1/2&    0&    0&    0&    0&    0 \\
1/2&    1&  1/2&    0&    0&    0&  1/2 \\
  0&  1/2&  1/2&  3/4&    0&    0&    0 \\
  0&    0&  3/4&  1/2&  3/4&    0&    0 \\
  0&    0&    0&  3/4&  1/2&  3/4&    0 \\
  0&    0&    0&    0&  3/4&  1/2&    0 \\
  0&  1/2&    0&    0&    0&    0&s_{16}\endbmatrix,\quad
s_{16}:=-\frac{3}{16}-
\frac{1}{8}\ooD d^2,
\]
where $s_2$ is as in~\eqref{eq.6A3.form}. As above,
$\discr_2\Fano_{\ttD}^d(\graph)$ depends on 
$\ooD\bmod4$
only. The minimal exceptional classes
are $\cl[8\bar\eta_{16}]=\cl[2\kk]$ and $\cl[\Ga_2]$ as in~\eqref{eq.6A3.classes} and,
as in the case of lines,
$\graph$ is not subgeometric unless $\ttD=4\bmod8$.
\endproof

\remark
In fact, in the \emph{non}-existence statements in this section, we do not
use the assumption $\ooD\gg d$. The reference to \autoref{lem.k.perp}
\emph{limits} the set of exceptional classes, \ie, for $\ooD$ small the
situation can only become \emph{worse}.
\endremark


\section{Type $\ttA{p-1}$ pencils for odd prime $p$}\label{S.A2}

We consider graphs of the type
$$
\graph := m\tA_{p-1} \oplus s\bA_1,
\quad p=3,5,\quad m\ge4,\quad s\le3.
$$
The vertices of the $i$-th $\tA_{p-1}$-subgraph of $\graph$
are denoted as $a_i^0,\dots,a_i^{p-1}$, where $i=1, \ldots, m$,
and the vertices of the $\bA_1$-components are
$b_1,\ldots,b_s$.
We put
$$
k:=\sum_{r=0}^{p-1}a_i^r=\const(i),\ \text{\cf. \eqref{eq.k}},  \qquad b:=b_1+ \dots +b_s
$$
and define
\begin{alignat*}3
\sa_i&:=a_i^1+a_i^2, &\quad \Ga_i&:=\tfrac13(a_i^1-a_i^2)
 &\quad \mbox{for $p=3$},\\
a_i^*&:=2a_i^1+3a_i^2+3a_i^3+2a_i^4, &  \Ga_i&:=\tfrac15(a_i^1+2a_i^2+3a_i^3+4a_i^4)
 &\quad \mbox{for $p=5$}.
\end{alignat*}
Note that $-a_i^*$ is the \emph{intrinsic polarization}
of the component in the sense of
\cite{degt:lines}.

\subsection{Preliminaries} \label{s.prelim}
Apart from
the lattice
$\FfNS$
generated by
\[*
h,k,a_i^r,b_j,\quad
i=1,\ldots,m,\
r=1,\ldots,p-1,\
j=1,\ldots,s,
\]
we
consider the auxiliary
lattice 
\[ \label{eq-barfano}
\barFfNS :=   (m\bA_{p-1} \oplus s\bA_1) \oplus  (\ZZ k + \ZZ \hbar),\quad
\mbox{where}\  \hbar:=h+d\sum_{i=1}^{m}\sa_i+\frac{d}2b  \, ,
\]
with the Gram matrix
\[*
m\bA_{p-1} \oplus s\bA_1 \oplus \bmatrix0&pd\\pd&\hbar^2\endbmatrix   \, .
\]
Obviously, for $s=0$ we have $\barFfNS = \FfNS$.
In general,
 \begin{itemize}
 	\item 
the discriminants of the two lattices
are of the same order, \ie,
 	 	\begin{equation} \label{eq-discrNS}
 	 	|\discr \FfNS| = |\discr \barFfNS| = p^{2+m} \cdot 2^s \cdot d^2
 	 	\end{equation}	
 	 	(indeed, this is an immediate consequence of the
definition of
$\hbar$);
 		\item for each prime $q > 2$ we have $\FfNS \otimes \Z_q  =  \barFfNS \otimes \Z_q$,
resulting in
 		\begin{equation} \label{eq-discriminants}
 		\discr_q \FfNS = \discr_q \barFfNS \, ;
  		\end{equation}  	
  		\item
one has
\[*
\OG_\hbar\bigl(\barFfNS\bigr)=\stab\hbar\subset\OG_h\bigl(\FfNS\bigr);
\]
\item
the orthogonal complements $k^\perp$ in the two lattices are equal; in
particular, we can use the concepts of shortest representatives
(\cf. \autoref{rem:smallperiod}) and exceptional classes/orbits (\cf.
\autoref{def.ex.orbit}).
 \end{itemize}
 The splitting \eqref{eq-barfano}
 results in the direct sum  decomposition
  \[ \label{eq-pdiscrdirsum}
  \discr_q\bigl(\FfNS\bigr) = \CD_q \oplus \CN_q \quad \mbox{ for } q > 2
  \]
  with
  ``the standard part'' $\CD_q:=\discr_q(m\bA_2 \oplus s\bA_1)$
  and  the group
  $$\CN_q := \discr_q(\ZZ k + \ZZ\hbar)$$
  that will play a crucial r\^{o}le in our
arguments.
For $q=p$
\eqref{eq-pdiscrdirsum} reads
 \[ \label{eq-discriminants}
 \discr_p\FfNS=\bigoplus_{i=1}^m\<\Ga_i\>\oplus\CN_p,
 \quad\Ga_i^2=\frac{p-1}{p}\bmod2\Z.
 \]

We start with a few useful observations. By the definition, we
have (still $q>2$)
\begin{equation} \label{eq-atmosttwo}
\ell(\CN_q)\le2.
\end{equation}
This group is nontrivial if and only if $\ppart{q}(pd)>0$,
and it is cyclic, \ie, $\ell(\CN_q)=1$, if and only if
$\ppart{q}(\hbar^2)=0$.

\lemma \label{lem-determinant-hkcyclic}
Let $q>2$ be a prime.
If  $\CN_q$  is cyclic, \ie,  $\ell(\CN_q)=1$,  then
$$
\det\CN_q = (-1) \cdot \hbar^2 \bmod(\Z_q\units)^2 \, .
$$
{In this case, all isotropic classes in $\CN_q$ are those of rational multiples of~$k$,
and they are all exceptional due to \autoref{lem.no.k}.}
\done
\endlemma

We need to make use of the discriminant class of $\hbar/q$, which is
well defined upon tensoring by $\Z_q$. However, to be on the safe
side, we introduce the vectors
\[*
\kk_q:=\frac{1}{q}k,\qquad \eta_q:=\frac{q+1}{q}\hbar\in\FfNS\dual;
\]
certainly, at this point we assume that $\ppart{q}(pd)>0$ and
$\ppart{q}(\hbar^2)>0$
{to make sure that the two vectors are indeed in the dual
lattice.}

\lemma \label{lem-determinant-hknoncyclic}
Let $q>2$ be a prime.
If $\ell(\CN_q)=2$, then
$$
\det\CN_q = (-1)  \bmod(\Z_q\units)^2  \, .
$$
In this case, the order~$q$ classes  in $\CN_q$   are
of the form
\[*
u\kk_q+v\eta_q,\quad u,v\in\Z/p.
\]
An isotropic class of this
form is exceptional if and only if $v=0$.
\done
\endlemma

In the sequel we will use the following general lemma on group forms.
\lemma\label{obs.K/K}
Let $\CK\subset\CM_q$ be an isotropic subgroup of order $q$
generated by a certain element $\Ga\in\CM_q$.
If there exists
another \emph{order~$q$} element $\Gb\in\CM_q$ such that
$\Ga\cdot\Gb\ne0\bmod\Z$,
then
\[*
\ell(\CK^\perp\!/\CK)=\ell(\CM_q)-2 \quad \mbox{and}\quad
\det (\CK^\perp\!/\CK)=-\det\CM _{q}   \bmod(\Z_q\units)^2.
\]
\endlemma

\proof
Abstractly, $\Ga$ and $\Gb$ generate $(\Z/q)^2$ with the Gram
matrix
\[*
\bmatrix0&\Ge/q\\\Ge/q&\ast  \endbmatrix,\quad\Ge\ne0\bmod q,
\]
of determinant $-1\bmod(\Z_q\units)^2$.
Since the form is obviously nondegenerate, $\Ga$ and $\Gb$ are independent in
$\CM_q$ and generate
an orthogonal direct
summand.
It is
this summand that disappears in
$\CK^\perp\!/\CK$.
\endproof

\observation \label{obs.cyclic}
Let $q=p$.
If $\CN_p$ is cyclic, as in \autoref{lem-determinant-hkcyclic},
the orbits of the action
of $O_h(\FfNS)$ on the isotropic
order~$p$ subgroups of
$\discr_p(\FfNS)$ are
\roster*
\item
the exceptional orbit
of $\langle\kappa\rangle$
and,
\item
typically, two
orbits of $\langle\Gapp_1\rangle,\langle\Gapp_2\rangle \subset \CD_p$,
one exceptional and
one not.
\endroster
(This fact is to be established on the case-by-case basis.
Note that
$\kappa+\Gapp_i\sim\Gapp_i$ and that any isotropic class of order
$p^r>p$ is exceptional, see \autoref{lem.no.k}.)
Hence, with one exception (see \autoref{s.8A2} below),
each $\graph$-admissible kernel  $\CK$ must be cyclic of order $p$.
Lemmata~\ref{obs.K/K} and~\ref{lem-determinant-hkcyclic}
imply that
\[* 
\ell\bigl(\discr_p\Fano_{\ttD}^d(\graph,\CK)\bigr) = \ell\bigl(\discr_p\FfNS\bigr) -2
\]
and $\discr_p\Fano_{\ttD}^d(\graph,\CK)$ satisfies
condition~\iref{Nikulin.p} of \autoref{th.Nikulin} if and only if
\[*
pm+s<21\quad\text{or}\quad
pm+s=21\ \,\text{and}\ \,\hbar^2=(-1)^m2^s\bmod p.
\]
\endobservation

\observation\label{obs.non-cyclic}
Still $q=p$.
If $\CN_p$ is not cyclic,
as in \autoref{lem-determinant-hknoncyclic},
in addition to
 the orbits of
  $\langle\kappa\rangle$, $\langle\Gapp_1\rangle$, $\langle\Gapp_2\rangle$
we also have that of $\<\eta_p+\Gb\>$
for some $\cl[\Gb]\in\CD_p$.
Since we assume $m\ge4$, we can redecompose
\[*
\CD_p=\underbrace{\biggl[\frac{2}p\biggr]\oplus
 \biggl[-\frac{2}p\biggr]}_{\CD'}\oplus
\underbrace{\biggl[\frac{2}p\biggr]\oplus
 \biggl[-\frac{2}p\biggr]}_{\CD''}\oplus\ldots
\]
so that $\cl[\Go_2]\in\CD'$. Thus, no matter what $\eta_p^2\bmod2\Z$ is, we can
find $\cl[\Gb],\cl[\Gg]\in\CD''$ (obviously orthogonal to $\cl[\Go_2]$)
with the property that
$(\eta_p+\Gb)^2=0\bmod2\Z$ and $(\eta_p+\Gb)\cdot\Gg\ne0\bmod\Z$.
Then, the kernel $\CK := \langle \cl[\Gapp_2],  \cl[\eta+\Gb]\rangle$ is
$\graph$-admissible and, applying \autoref{obs.K/K} twice
and \autoref{lem-determinant-hknoncyclic}, we find that
$$
\ell\bigl(\discr_p\Fano_{\ttD}^d(\graph,\CK)\bigr) = \ell\bigl(\discr_p\FfNS\bigr) - 4
$$
and a \emph{sufficient} condition for
$\discr_p\Fano_{\ttD}^d(\graph,\CK)$ to satisfy
\autoref{th.Nikulin}\iref{Nikulin.p} is
\[*
pm+s<22\quad\text{or}\quad
pm+s=22\ \,\text{and}\ \,(-1)^{m+1}=2^s\bmod p.
\]
The last relation holds in all applications below;
therefore, we do not need to engage into a study of isotropic elements of
higher order (\cf. \autoref{lem.q}).
\endobservation

\subsection{Other torsion in the discriminant}\label{s.other.torsion}
A necessary
realizability condition that cannot be achieved by a
finite index extension of the lattice is \autoref{th.Nikulin}\iref{Nikulin.r}.
Therefore, we
always assume that
\[*
(p-1)m+s\le18.
\]

If $q\ne p$ is an odd prime, then $\discr_q\FfNS=\CN_q$ and
$\ell(\CN_q)\le2$, see~\eqref{eq-atmosttwo}.
We need the following arithmetic condition: given a prime~$q$,
a pair $(h,d)$ of rational numbers is called \emph{$q$-balanced} if either
\[*
\ppart{q}(h)\ge\ppart{q}(d)\quad\text{or}\quad \text{$\ppart{q}(h)$ is even},
\]
\cf. the statement of \autoref{thm.main}.

\lemma\label{lem.q}
Let $\ooD\gg d>0$ and $q\ne p$ a prime. Then
$\graph$ has a $(\bmod\,q)$-geometric kernel if and only if either
\roster*
\item
$q=1\bmod4$, or
\item
$(p-1)m+s\le17$, or
\item
the pair $(\hbar^2,d)$ is $q$-balanced.
\endroster
\endlemma

\proof
The case $(p-1)m+s\le17$ is trivial (as then we need $\ell(\CN_q)\le3$), and
so is the case $q=1\bmod4$ (as then $-1\bmod q$ is a quadratic residue).
Therefore, in the rest of the proof we assume  $(p-1)m+s=18$ and $q=3\bmod4$.
Condition~\iref{Nikulin.p} of \autoref{th.Nikulin} is
immediately satisfied whenever
\roster*
\item
$\CN_q=0$, \ie, $\ppart{q}(d)=0$, or
\item
$\ell(\CN_q)=1$, \ie, $\ppart{q}(\hbar^2)=0<\ppart{q}(d)$.
\endroster
Note that, in the last two cases, the pair $(\hbar^2,d)$ is $q$-balanced.

Thus, assume that $u:=\ppart{q}(d)>0$ and $v:=\ppart{q}(\hbar^2)>0$ and let
$t:=\min\{u,v\}$. Then,
\[*
\CN_q=\biggl[-\frac{2\Ge}{q^{2u-t}}\biggr]\oplus\biggl[\frac{2\Ge}{q^t}\biggr],
\]
generated by some classes $x,y$. If $t=u$, this group has at least two
cyclic isotropic subgroups~$\CK$ of order~$u$, \viz. those generated by
$x\pm y$, and at least one of them does not contain $\kk_q$ and, hence, is
$\graph$-admissible. In this case, $\CK^\perp\!/\CK=0$.

If $t=v<u$ is even, $t=2r$, we can take $y=\cl[\hbar/q^t]$ to see that the
isotropic subgroup~$\CK$ generated by $q^ry$ is $\graph$-admissible and
$\ell(\CK^\perp\!/\CK)=1$.

Finally, assume that $t=v<u$ is odd, $t=2r+1$, \ie, the only case where the
pair $(\hbar^2,d)$ is not $q$-balanced. The subgroup
\[*
\CK=\CN_q^r:=\bigl\{z\in\CN_q\bigm|q^rz=0\bigr\}
\]
is
isotropic, but it does not suffice: we have
$\CK^\perp\!/\CK\same(\Z/q^{2u-2r})\oplus(\Z/q)$. Hence, we need at least
one isotropic class~$\Gg$ of order $q^{r+1}$. Modulo $\CN_q^r$, all such classes
are multiples of $x$; since, on the other hand, all rational multiples
of~$k$ are isotropic, any ``suitable'' kernel
must contain the exceptional class $\cl[\kk_q]=\cl[q^rx]$.
\endproof

\observation\label{obs.q=2}
We can also compute the $2$-torsion. If $d$ is odd, then
\[*
\discr_2\FfNS=s\biggl[-\frac{1}2\biggr],\quad\text{generated by}\ \,
 \Gb_i:=\frac12(k+b_i)
\]
by~\eqref{eq-discrNS} and \autoref{rem:compdiscr}. It is odd and, since
$s\le3$, has no isotropic subgroups. Thus, condition~\iref{Nikulin.2} in
\autoref{th.Nikulin} is satisfied if and only if
\[*
(p-1)m+2s\le20.
\]
If $d$ is even, then \eqref{eq-barfano} is a valid change of basis
\emph{over~$\Z$} and we can use~\eqref{eq-pdiscrdirsum} for $q=2$. We
conclude that $\discr_2\FfNS$ is of length $s+2$ and odd (unless $s=0$ and
$\ppart2(\hbar^2)\ge2$) and the inequality
\[*
(p-1)m+2s\le18
\]
is sufficient for condition~\iref{Nikulin.2} in
\autoref{th.Nikulin}. For some values of $\ooD$, $d$ it may not be necessary;
the tedious analysis of isotropic subgroups is left to the reader, \cf. the
proof of \autoref{lem.q}.
\endobservation

%


\subsection{The pencil $8\ttA2$}\label{s.8A2}
For $\graph := 8\tA_2$
we have
\[*
\hbar^2=2n+16d^2,
\]
see \eqref{eq-barfano}.
The following claim strengthens
\cite[Proposition~6.1]{rams.schuett.smoothodd};
at the same time, it serves as an example of the machinery explained in
this section.

\proposition\label{th.8A2}
For integers $\ooD\gg d>1$, the statements below are equivalent\rom:
\roster*
\item
the graph $\graph:=8\tA_2$ is geometric in degree~$\ttD$ with respect
to~$d$\rom;
\item
the discriminant group $\CN_3:=\discr_3(\Z k+\Z\hbar)$ is \emph{not} cyclic\rom;
\item
$\hbar^2=0\bmod3$\rom;
\item
$\ttD=2d^2\bmod3$.
\endroster
\endproposition

\proof
The last three statements are obviously equivalent.

We have $\rank\FfNS=18$ and
$\ell\bigl(\discr_3\FfNS\bigr)=8+\ell(\CN_3)\ge9$.

If $\CN_3$ is cyclic, \ie, $\hbar^2\ne0\bmod3$, \cf. \autoref{obs.cyclic},
then, apart from $\<\kk\>$,
there are two
$\OG_h\bigl(\FfNS\bigr)$-orbits of isotropic order-$3$ subgroups, with
representatives generated by the classes of
\[
\Go_1:=\Ga_1+\Ga_2+\Ga_3\quad\text{or}\quad
\Go_2:=\Ga_1+\Ga_2+\Ga_3+\Ga_4+\Ga_5+\Ga_6.
\label{eq.8A2}
\]
(We freely use the notation of \autoref{s.prelim}.)
The former is exceptional, the latter is not, as its shortest representatives
have square~$(-4)$.
This is the only case where we can have a $\graph$-admissible kernel~$\CK$ of
length~$2$, \eg, the one generated by
\[*
\cl[\Go_2]\quad\text{and}\quad
\cl[\Go_2']\ni\Ga_3+\Ga_4-\Ga_5-\Ga_6+\Ga_7+\Ga_8.
\]
Hence, we can reduce $\ell\bigl(\discr_3\FfNS\bigr)=9$ down to~$5$, which is
still not enough for condition~\iref{Nikulin.p} in \autoref{th.Nikulin}.

If $\CN_3$ is not cyclic, $\hbar^2=0\bmod3$, then, arguing as in
\autoref{obs.non-cyclic}, we can find a $\graph$-admissible kernel~$\CK$ of
length~$3$, \viz. the one generated by $\cl[\Go_2]$, $\cl[\Go_2']$, and an element
of the form $\cl[\eta_3+\Gb]$, $\Gb\subset\CD_3$, and reduce
$\ell\bigl(\discr_3\FfNS\bigr)=10$ by six units, down to~$6$. Since
\autoref{obs.K/K} is applied \emph{three} times,
\autoref{lem-determinant-hknoncyclic} yields
\[
\det\nolimits_3(\CK^\perp\!/\CK)=1\bmod(\Z_3\units)^2,
\]
which fits condition~\iref{Nikulin.p} in \autoref{th.Nikulin}.

For the other primes $q\ne3$, see \autoref{lem.q} and \autoref{obs.q=2}.
\endproof

\remark
We do not
assert
the uniqueness of
a family
of \KK-surfaces
that contain an $8\tA_2$ configuration
we just constructed.
The primes $q\mathrel|d$, $q\ne3$,
often result in several admissible extensions.

The assumption $h^2\gg0$ may be relevant as occasionally
there are $(-2)$-vectors orthogonal to~$h$ but not to~$k$.
\endremark

\subsection{Other type~$\ttA2$ pencils}\label{s.A2}
The other three type~$\tA_2$ pencils that we need for the proof of
\autoref{thm.main} are straightforward:
for $\discr_q\FfNS$, we use
\roster*
\item
the identity $\discr_q\FfNS=\CN_q$ and \autoref{lem.q} for $q>3$,
\item
Observations~\ref{obs.cyclic},~\ref{obs.non-cyclic} and
the list~\eqref{eq.8A2} for $q=3$, and
\item
\autoref{obs.q=2} for $q=2$.
\endroster
For the readers convenience, we compute
\[*
\hbar^2=2n+\left(2m+\frac{s}{2}\right)d^2\quad\text{for}\ \,
\graph:=m\tA_2+s\bA_1,
\]
whereupon we merely state the ultimate results.

\proposition\label{th.7A2}
For integers $\ooD\gg d>1$, the statements below are
equivalent\rom:
\roster*
\item
the graph $\graph:=7\tA_2$ is geometric in degree~$\ttD$ with respect
to~$d$\rom;
\item
$\hbar^2\ne2\bmod3$\rom;
\item
$2n\ne d^2-1\bmod3$.
\qed
\endroster
\endproposition

\proposition\label{th.7A2A1}
For integers $\ooD\gg d>1$, the statements below are
equivalent\rom:
\roster*
\item
the graph $\graph:=7\tA_2\oplus\bA_1$ is geometric in degree~$\ttD$ with respect
to~$d$\rom;
\item
the discriminant group $\CN_3:=\discr_3(\Z k+\Z\hbar)$ is \emph{not} cyclic\rom;
\item
$\hbar^2 = 0\bmod3$\rom;
\item
$2n = 2d^2\bmod3$.
\qed
\endroster
\endproposition

\proposition\label{th.6A2+3A1}
For integers $\ooD\gg d>1$, the statements below are
equivalent\rom:
\roster*
\item
the graph $\graph:=6\tA_2\oplus3\bA_1$ is geometric in degree~$\ttD$ with respect
to~$d$\rom;
\item
$\hbar^2\ne 1 \bmod3$\rom;
\item
$2n\ne1\bmod3$.
\qed
\endroster
\endproposition

\subsection{The pencil $4\ttA4+2\bA\sb1$}\label{S.4A4+2A1}

In this section, we do have to \emph{assume that $d$ is odd,
\cf. \autoref{obs.q=2}.} We keep using the notation of \autoref{s.prelim}.


\proposition\label{th.4A4+2A1}
Let $\ooD\gg d>1$ be integers with $d$ \emph{odd}. Then, the following
statements are equivalent\rom:
\roster*
\item
the graph $\graph:=4\tA_4\oplus2\bA_1$ is geometric in degree~$\ttD$ with respect
to~$d$\rom;
\item
the discriminant group $\CN_5:=\discr_5(\Z k+\Z\hbar)$ is \emph{not} cyclic
and the pair $(\hbar^2,d)$ is $q$-balanced for each prime $q=3\bmod4$\rom;
\item
$\hbar^2=0\bmod5$ and $(\hbar^2,d)$ is $q$-balanced for each prime $q=3\bmod4$\rom;
\item
$2n=4d^2\bmod5$ and $(\ttD,d)$ is $q$-balanced for each prime $q=3\bmod4$.
\endroster
\endproposition

\proof
We have
\[*
\hbar^2 = \ttD + 41 d^2,
\]
and it is clear that the last three statements are equivalent.

The $5$-torsion is handled by Observations~\ref{obs.cyclic}
and~\ref{obs.non-cyclic}: shortest representatives in $\cl[\pm\Ga_i]$ and
$\cl[\pm2\Ga_i]$ have squares $-4/5$ and $-6/5$, respectively, and the two
orbits in \autoref{obs.cyclic} are those of
\[*
\Go_1:=\Ga_1+2\Ga_2,\qquad \Go_2:=\Ga_1+2\Ga_2+\Ga_3+2\Ga_4;
\]
the former is exceptional, the latter is not. We conclude that a
$(\bmod\,5)$-geometric kernel exists if and
only if $\CN_5$ is not cyclic.

For the other primes $q\ne5$, we use \autoref{lem.q} and \autoref{obs.q=2}:
since $\rank\FfNS=20$ is maximal, the $q$-balancing condition
does play a r\^{o}le.
\endproof

\example\label{ex.low.degree}
Let $d=3$. Consider the pencil $\graph:=4\tA_4\oplus2\bA_1$ and the lattice
$N:=\Fano_\ttD^3(\graph,\CK)$, where $\CK$ is a geometric kernel containing
$\Gl:=k/3$, see~\eqref{eq.k}. The latter vector is represented by the roots
of the form
\[*
l_i^r:=-\Gl+a_i^{2r}+a_i^{2r+1}\in\root_1(N,h),\quad i=1,\ldots,4,\ r=0,\ldots,4.
\]
(Here and below, the upper index should be understood $\bmod\,5$.
We assume that $\ooD\gg0$ and, hence, there are no other
vectors of interest.) Thus, we do obtain a $\ttD$-polarized
$K3$-surface~$X$ with an elliptic pencil with four type~$\tA_4$ fibers, but
the fiber components are the lines $l_i^r$, whereas the anticipated
original cubics are split:
\[*
a_i^r=l_i^{3r}+l_i^{3r+1}+l_i^{3r+2}.
\]
Of course, the two type~$\bA_1$ fibers remain cubics.

In this particular case, there is another, geometric kernel; we chose this
example because it is more transparent geometrically. To illustrate the
importance of the $q$-balancing condition, one would have to take, \eg, $d=9$
and $2n=24\bmod90$: instead of degree~$9$ curves, each fiber would
\emph{have} to consist of
five twisted cubics.
\endexample

\section{Proof of \autoref{thm.main} and \autoref{cor-miyaoka-classification}} \label{S.proof}

\subsection{Proof of \autoref{thm.main}}\label{proof.main}
Let
$\ooD\gg d$.
By \cite[Proposition A]{Miyaoka:lines}, the maximal number of smooth
rational degree~$d$ curves on a \smooth
$\ttD$-polarized $K3$-surface
does not exceed $24$.
If either $d$ or $\ooD$ is even, the existence of
surfaces with $24$ smooth
rational degree~$d$ curves follows from \cite[Theorem~1.1]{rams.schuett.smoothodd}.
Thus, we can assume that both $\ooD$ and $d$ are odd.

For all $\ooD\gg d>0$,
Propositions~\ref{th.7A2} and~\ref{th.6A2+3A1} give
us surfaces with $21$ curves; therefore, we concentrate on those with at
least $22$ curves,
in which case (still assuming that $\ooD\gg d$),
all smooth rational curves of degree \emph{at most~$d$} must be fiber
components of an elliptic pencil, see \cite[Theorem~1.1(ii)]{rams.schuett:24}
or~\cite{degt:lines}.
In other words,
to construct a \smooth $K3$-surface with $N=22$, $23$, or~$24$
smooth rational curves, we need a graph $\tilde{\graph}$ appearing on
Shimada's list~\cite{Shimada:ellipticK3} of elliptic pencils, with a certain
number $N\le F\le24$ vertices, and such that, upon assigning color~$d$ to $N$
vertices and a color \emph{other than~$d$} to the remaining $\Gd:=F-N$
vertices, the \emph{weighted degree} $\sum_{v\in\Sigma}n_v\deg(v)$,
\cf.~\eqref{eq.k}, of the components $\Sigma\subset\tilde{\graph}$ is
constant.
Then our graph~$\graph$ is obtained from~\smash{$\tilde{\graph}$} by removing all
vertices of degree other than~$d$.

Denote by~$f$ the ``expected'' plain degree
$\deg(\graph)$, \ie, the common degree of the components of~\smash{$\tilde{\graph}$} that
are to remain parabolic.

If $\Gd=0$, then $\deg(\Sigma)=f$ for \emph{each} component
$\Sigma\subset\tilde{\graph}$, leaving but four pencils:
\[ \label{eq-fourpencils}
\graph=\tilde{\graph}=12\tA_1,\quad 11\tA_1,\quad 8\tA_2,\quad\text{or}\quad 6\tA_3,
\]
see Propositions~\ref{lem:A1}, \ref{th.8A2}, and~\ref{lemma.6tA3}.

If $\Gd=1$, then $\tilde{\graph}$ has one distinguished component $\Sigma'$
of $\deg(\Sigma')<d$ (as otherwise the vertex removed would have to have the
same degree~$d$), and then $\graph$ is obtained from~\smash{$\tilde{\graph}$}
by ``stripping'' the $\tilde{\ }$ from~$\Sigma'$:
\[*
\graph=7\tA_2\oplus\bA_1\quad\text{or}\quad5\tA_3\oplus\bA_2,
\]
see Propositions~\ref{th.7A2A1} and~\ref{lemma.5tA3plusbA2}.

Finally, if $\Gd=2$, then either $\tilde{\graph}$ has two components
$\Sigma',\Sigma''$ of $\deg(\Sigma^*)<d$:
\[*
\graph=5\tA_3\oplus2\bA_1\quad\text{or}\quad4\tA_4\oplus2\bA_1,
\]
see Propositions~\ref{lemma.5tA3plus2bA1} and~\ref{th.4A4+2A1},
or $\tilde{\graph}$ has one distinguished component~$\Sigma'$, of any degree
(assuming $d>1$),
with $3\le\ls|\Sigma'|\le f+1$:
\[
\graph=7\tA_2\oplus(\tA_2\mapsto\bA_1)\quad\text{or}\quad
5\tA_3\oplus(\tA_3\mapsto\text{$\bA_2$ or $2\bA_1$}).
\label{eq.not-spanned}
\]
It is this last case where we would have to use a version of
\autoref{thm:K3} that does not assume that
$\NS(X)$ is spanned by degree~$d$ curves (\cf. \autoref{rem.Schur}). Fortunately, all three
configurations have
{just been discussed as obtained by other means,
satisfying the assumptions of \autoref{thm:K3}.}

Summarizing, \autoref{th.8A2} gives us an extra family (in addition to the
case $\ooD d=0\bmod2$, see \cite[Theorem~1.1]{rams.schuett.smoothodd})
of pairs $(\ooD,d)$
realized by a \smooth $K3$-surface with $24$ curves, whereupon
Propositions~\ref{lem:A1}, \ref{lemma.5tA3plus2bA1}, \ref{lemma.5tA3plusbA2},
and~\ref{th.7A2A1} assert that the other type $\tA_p$ graphs~$\graph$,
$p=2,3,4$, with $\ls|\graph|\ge22$
do not give us anything new. The only case where $22$, but not more, curves
can be obtained is $4\tA_4\oplus2\bA_1$ that is considered in
\autoref{th.4A4+2A1}.
\qed
\subsection{Proof of  \autoref{cor-miyaoka-classification}}\label{proof.cor-miyaoka-classification}
{Assume that a complex \KK-surface contains}
$24$ rational degree~$d$ curves
{and $2n\gg d$}.  By \cite[Theorem~1.1(ii)]{rams.schuett:24}
(see also \cite{degt:lines})  the rational curves are  components of singular
fibers of a genus one fibration on $X$.
Since
{the topological} Euler--Poincar\'{e} characteristic of $X$ is $24$ and general
fiber of the fibration
is smooth, \emph{all} components of singular fibers must be irreducible rational
degree~$d$ curves (see \eg~\cite[(5.4)]{rams.schuett:24}).

Thus, if one singular fiber is irreducible, then
{so are the others}
(all fibers are of degree $d$). This results in the configuration
$24\tA_0^*$.
{This configuration exists for all pairs $(n,d)$
provided that $d>2$, see
\cite[\S\,10]{rams.schuett:24}.}

Otherwise, all singular fibers are reducible
{and, hence, all}
$24$ curves are smooth.
{In view of \eqref{eq-fourpencils},}
Propositions~\ref{lem:A1}, \ref{th.8A2}, and~\ref{lemma.6tA3}
complete the proof.
\qed

\remark[\cf.~ \autoref{cor-miyaoka-classification}] \label{rem:even.nd}
Let $\ooD \gg d$ with $nd$ \emph{even}.
For
completeness,
we collect
all known {(to us)} facts about
{the} existence of configurations of $24$ irreducible rational degree~$d$
curves on a complex \KK-surface of degree $2\ooD$ in this case.
\roster*
\item $24\tA_0^*$ exists for all $(n,d)$ provided {that} $d>2$,
see \cite[\S\,10]{rams.schuett:24};
\item $12\tA_1$ exists for all $(n,d)$
by \cite[Theorem~1.1]{rams.schuett.smoothodd};
\item
{$8\tA_2$ exists
if and only if $\ooD=d^2\bmod3$ by \autoref{th.8A2};}
\item
{we conjecture that $6\tA_3$  exists if and only if $\ooD = 2d^2 \bmod 4d$,
\cf.~\cite[\S\,11.1]{rams.schuett:24}.}
\endroster
{The proof of the last conjecture, which seems too technical for the
present paper, will appear elsewhere.}
\endremark

{
\let\.\DOTaccent
\def\cprime{$'$}
\bibliographystyle{amsplain}

}
\end{document}